\documentclass[11pt,letterpaper,reqno]{amsart}

\usepackage{amssymb}
\usepackage{amsmath}
\usepackage{amsthm}
\usepackage{amsfonts}
\usepackage{bbm}

\usepackage{bookmark}
\usepackage{hyperref}
\hypersetup{pdfstartview={FitH}}

\newtheorem{theorem}{Theorem}

\numberwithin{equation}{section}

\begin{document}

\title[Popular differences for right isosceles triangles]{Popular differences for right isosceles triangles}

\author[Vjekoslav Kova\v{c}]{Vjekoslav Kova\v{c}}
\address{Vjekoslav Kova\v{c}, Department of Mathematics, Faculty of Science, University of Zagreb, Bijeni\v{c}ka cesta 30, 10000 Zagreb, Croatia}
\email{vjekovac@math.hr}


\subjclass[2020]{
Primary
11B30; 
Secondary
05D05} 

\keywords{Roth's theorem, Szemer\'{e}di's theorem, Fourier transform, polynomial method}

\begin{abstract}
For a subset $A$ of $\{1,2,\ldots,N\}^2$ of size $\alpha N^2$ we show existence of $(m,n)\neq(0,0)$ such that the set $A$ contains at least $(\alpha^3 - o(1))N^2$ triples of points of the form $(a,b)$, $(a+m,b+n)$, $(a-n,b+m)$.
This answers a question by Ackelsberg, Bergelson, and Best.
The same approach also establishes the corresponding result for compact abelian groups.
Furthermore, for a finite field $\mathbb{F}_q$ we comment on exponential smallness of subsets of $(\mathbb{F}_q^n)^2$ that avoid the aforementioned configuration.
The proofs are minor modifications of the existing proofs regarding three-term arithmetic progressions.
\end{abstract}

\maketitle


\section{Introduction}
The main result of this note deals with a configuration consisting of three distinct points in $\mathbb{Z}^2$ of the form
\begin{equation}\label{eq:triangle}
(a,b), \ (a+m,b+n), \ (a-n,b+m).
\end{equation}
In other words, we are interested in triples of points forming vertices of a right isosceles triangle.

\begin{theorem}\label{thm:main}
For every $\varepsilon\in(0,1]$ there exists a positive integer $N_0(\varepsilon)$ such that the following holds: for every integer $N\geq N_0(\varepsilon)$ and a set $A\subseteq\{1,2,\ldots,N\}^2$ with density $\alpha:=|A|/ N^2$ there exists a pair $(m,n)\in\mathbb{Z}^2\setminus\{(0,0)\}$ such that
\[ \big| A \cap \big(A - (m,n)\big) \cap \big(A - (-n,m)\big) \big| \geq (\alpha^3 - \varepsilon)N^2, \]
i.e., $A$ contains at least $(\alpha^3 - \varepsilon)N^2$ triples of points of the form \eqref{eq:triangle} for the same fixed value of $(m,n)\neq(0,0)$.
\end{theorem}

Theorem~\ref{thm:main} confirms the conjecture of Ackelsberg, Bergelson, and Best formulated in a preprint version of the paper \cite{ABB21:tri}, later acknowledged as \cite[Theorem~1.19]{ABB21:tri}. In fact, these authors have established a certain infinitary result dealing with the same point configuration \eqref{eq:triangle} using multiple recurrence techniques from ergodic theory, while our approach is combinatorial. 
However, note that the Furstenberg correspondence principle does not apply to this type of problems, so Theorem~\ref{thm:main} cannot be automatically deduced from \cite{ABB21:tri}.

Analogous results for three-term and four-term arithmetic progressions,
\begin{equation}\label{eq:ap}
a, \ a+d, \ a+2d, \ ( a+3d, )
\end{equation}
had been conjectured by Bergelson, Host, and Kra \cite{BHK05:34ap}, and they were established by Green \cite{G05:3ap} and Green and Tao \cite{GT10:4ap}, respectively. These results are often called ``popular difference'' strengthenings of the theorems of Roth \cite{R53:Roth} and Szemer\'{e}di \cite{S69:Szem}. They find a nonzero $d$ such that the number of progressions \eqref{eq:ap} in the set $A\subseteq\{1,2,\ldots,N\}$ with difference $d$ comes close to the theoretical maximum, obtained by constructing $A$ at random.
In the appendix of \cite{BHK05:34ap}, Ruzsa gave a construction showing that arithmetic progressions of length five or more no longer share this property.
Moreover, Mandache \cite{M18:cor} showed that the corresponding result also fails for corners, i.e., patterns of the form \eqref{eq:triangle} with $n=0$, while Fox, Sah, Sawhney, Stoner, and Zhao \cite{FSSSZ20:cor} found the sharp asymptotics for this pattern.
Configurations \eqref{eq:triangle} allow an additional degree of freedom, as right isosceles triangles can now be rotated.
Let us also remark that Sah, Sawhney, and Zhao \cite{SSZ20:pat} characterized all point configurations for which a popular difference result holds, but they only allowed dilations and translations in pattern formation, so their discussion does not include rotated right isosceles triangles.

We will establish a more general result in the context of $2$-divisible compact Hausdorff abelian groups $(\mathbb{G},+)$, in perfect analogy with \cite[Theorem~3]{T14:Roth} from an expository note by Tao. Recall that a group $\mathbb{G}$ is \emph{$2$-divisible} if the map $\mathbb{G}\to\mathbb{G}$, $x\mapsto 2x$ is surjective. This is the case with finite cyclic groups of odd order, which will be sufficient for the intended application to Theorem~\ref{thm:main}.
The \emph{Bohr set} $B^{\mathbb{G}}(S,\rho)$ associated with a finite set of frequencies $S\subseteq\widehat{\mathbb{G}}$ and a radius $\rho\in(0,1]$ is
\[ B^{\mathbb{G}}(S,\rho) := \Big\{x\in\mathbb{G} \,:\, \max_{\xi\in S} \|\xi(x)\|_{\mathbb{R}/\mathbb{Z}}<\rho \Big\}. \]
The group $\mathbb{G}$ is equipped with its probability Haar measure $\mu$, which leads us to the algebra $\textup{L}^1(\mathbb{G})$ with respect to convolution.
We also define the \emph{tent function} $\nu_{S,\rho}$ over $B^{\mathbb{G}}(S,\rho)$ as a constant multiple of
\[ x \mapsto \Big(1 - \rho^{-1}\max_{\xi\in S}\|\xi(x)\|_{\mathbb{R}/\mathbb{Z}}\Big)_{+} \]
normalized so that $\|\nu_{S,\rho}\|_{\textup{L}^1(\mathbb{G})}=1$ and we set
$\chi_{S,\rho} := \nu_{S,\rho} \ast \nu_{S,\rho}$.
Finally, for any three functions $f,g,h\in\textup{L}^{\infty}(\mathbb{G}^2)$ and a ``weight'' function $\chi\in\textup{L}^{\infty}(\mathbb{G})$ define the \emph{weighted counting form}
\begin{equation}\label{eq:defLambda}
\Lambda_{\chi}(f,g,h) := \int_{\mathbb{G}} \int_{\mathbb{G}} \int_{\mathbb{G}} \int_{\mathbb{G}} f(x,y) \,g(x+s,y+t) \,h(x-t,y+s) \,\chi(s) \,\chi(t) \,\textup{d}\mu(s) \,\textup{d}\mu(t) \,\textup{d}\mu(x) \,\textup{d}\mu(y).
\end{equation}

\begin{theorem}\label{thm:compact}
For every $\varepsilon\in(0,1]$ there exists $R(\varepsilon)\in(0,\infty)$ such that the following holds: for every $2$-divisible compact Hausdorff abelian group $(\mathbb{G},+)$ with the Haar probability measure $\mu$ and a Borel-measurable function $f\colon\mathbb{G}^2\to[0,1]$ there exist a set of frequencies $S\subseteq\widehat{\mathbb{G}}$ with cardinality $|S|\leq R(\varepsilon)$ and a radius $\rho\in[R(\varepsilon)^{-1},\varepsilon]$ such that
\[ \Lambda_{\chi_{S,\rho}}(f,f,f) \geq \Big(\int_{\mathbb{G}}\int_{\mathbb{G}} f(x,y)\,\textup{d}\mu(x)\,\textup{d}\mu(y)\Big)^3 - \varepsilon .\]
In the case of $\mathbb{G}=\mathbb{Z}/N\mathbb{Z}$ we can choose $S,\rho$ so that $\chi_{S,\rho}$ is supported in
$\{-2\lfloor \rho N\rfloor,\ldots,2\lfloor \rho N\rfloor\} + N\mathbb{Z}$.
\end{theorem}

Derivation of Theorem~\ref{thm:main} from Theorem~\ref{thm:compact} is standard by embedding $\{1,2,\ldots,N\}$ in a slightly larger cyclic group of odd order and taking $f=\mathbbm{1}_A$. We have an upper bound on the height of the function $\chi_{S,\rho}$ in question (by the standard ``volume packing bound'' for Bohr sets \cite{T14:Roth}), so the trivial ``difference'' $(s,t)=(0,0)$ contributes negligibly to the count. The additional support requirement for $\chi_{S,\rho}$ makes the impact of wraparounds negligible too.

Let us remark that a variant of Theorem~\ref{thm:compact} in a more general context of matrix patterns was established by Berger, Sah, Sawhney, and Tidor \cite{BSST21:matrix} shortly after the appearance of this note. This resulted in an independent resolution to the aforementioned question of Ackelsberg, Bergelson, and Best.

From the previous discussion one gets an impression that configuration \eqref{eq:triangle} is easier than the previously mentioned arithmetic configurations. In this paper we argue that it is certainly not more difficult than the three-term arithmetic progression and we treat it in close analogy with it. This similarity in the ``Fourier complexity'' has also been noted by several other authors \cite{B14:PhD,BSST21:matrix,SS21:trian}.

Now let $\mathbb{G}$ be a finite abelian group with $N$ elements.
If one asks about the size of the largest subset of $\{1,2,\ldots,N\}^2$ or $\mathbb{G}^2$ that does not contain the configuration \eqref{eq:triangle}, then a bound $N^2/(\log\log N)^c$ comes from Shkredov's results \cite{Sh06:cor,Sh09:cor} on corners, i.e., by considering $n=0$ only. A stronger result, for rotated squares, has been established by Prendiville \cite{P15:matrix}.
An improvement is possible in the form of a polylogarithmic bound $N^2/(\log N)^c$ for sets avoiding \eqref{eq:triangle}. A standard density argument can be based on a ``uniformity'' estimate that follows from formula \eqref{eq:Four} below:
\[ |\Lambda_{\mathbbm{1}}(f,g,h)| \leq \min\big\{ \big\|\widehat{f}\big\|_{\ell^\infty(\widehat{\mathbb{G}}^2)}, \big\|\widehat{g}\big\|_{\ell^\infty(\widehat{\mathbb{G}}^2)}, \big\|\widehat{h}\big\|_{\ell^\infty(\widehat{\mathbb{G}}^2)} \big\} \]
for Borel-measurable functions $f,g,h\colon\mathbb{G}^2\to[-1,1]$, at least when $N=|\mathbb{G}|$ is odd.
We do not discuss this argument in detail, as, in fact, Bloom \cite[Theorem~2.21]{B14:PhD} showed a more general result, while studying certain multidimensional patterns.
Finally, Shkredov and Solymosi \cite{SS21:trian} recently posed a challenge of breaking the logarithmic barrier in this polylogarithmic bound for rotated triangles in $\{1,2,\ldots,N\}^2$; see \cite[Conjecture~2]{SS21:trian} and \cite[Remark~14]{SS21:trian}. 
The problem can be reformulated in terms of either three-term matrix progressions,
\[ \begin{bmatrix} a \\ b \end{bmatrix},\quad 
\begin{bmatrix} a \\ b \end{bmatrix} + \begin{bmatrix} 1 & 0 \\ 0 & 1 \end{bmatrix} \begin{bmatrix} m \\ n \end{bmatrix},\quad
\begin{bmatrix} a \\ b \end{bmatrix} + \begin{bmatrix} 0 & -1 \\ 1 & 0 \end{bmatrix} \begin{bmatrix} m \\ n \end{bmatrix}, \]
see \cite{B14:PhD,BSST21:matrix}, or three-term progressions in the Gaussian integers $\mathbb{Z}[\mathbbm{i}]$, 
\[ z, \quad z+w, \quad z+\mathbbm{i}w, \]
but the literature on both seems to be scarce.
In fact, the author is not aware of any asymptotically nontrivial lower bounds. Maximal sizes of sets in the grid $\{1,2,\ldots,N\}^2$ that do not contain \eqref{eq:triangle} are known for $N\leq 11$ and are listed as sequence A271906 in the encyclopedia OEIS \cite{OEIS}.

Instead, we will comment on the finite field case, i.e., when $\mathbb{G}=\mathbb{F}_q^n$, where $\mathbb{F}_q$ denotes the finite field with $q$ elements.

\begin{theorem}\label{thm:expsmall}
For every prime power $q$ there exists a number $c_q\in(0,q)$ such that the following holds for every positive integer $n$: if a set $A\subseteq(\mathbb{F}_q^n)^2$ does not contain a triple of distinct points \eqref{eq:triangle} with $a,b,m,n\in\mathbb{F}_q^n$, then its cardinality needs to satisfy the bound $|A|\leq 3c_q^{2n}$.
\end{theorem}

We will see that the polynomial method by Croot, Lev, and Pach \cite{CLP17:cap} and Ellenberg and Gijswijt \cite{EG17:cap} devised for three-term arithmetic progressions works here without any major modifications. In fact, we will follow its elegant reformulation in terms of the tensor slice rank due to Tao \cite{T16:poly}. This points out yet another similarity between configurations \eqref{eq:triangle} and \eqref{eq:ap}.

The proof of Theorem~\ref{thm:compact} is given in Section~\ref{sec:regularity}, while the proof of Theorem~\ref{thm:expsmall} is given in Section~\ref{sec:expsmall}. As we have already said, both of them will be minor modifications of the existing proofs.


\section{Proof of Theorem~\ref{thm:compact}}
\label{sec:regularity}
Before we begin with the proof, let us observe formulae
{\allowdisplaybreaks
\begin{subequations}
\begin{align}
\Lambda_{\chi}(f,g,h) & = \sum_{\xi,\zeta,\xi',\zeta'\in\widehat{\mathbb{G}}} \widehat{f}(-\xi-\xi',-\zeta-\zeta') \,\widehat{g}(\xi,\zeta) \,\widehat{h}(\xi',\zeta') \,\widehat{\chi}(-\xi-\zeta') \,\widehat{\chi}(-\zeta+\xi'), \label{eq:weiFour1} \\
\Lambda_{\chi}(f,g,h) & = \sum_{\xi,\zeta,\xi',\zeta'\in\widehat{\mathbb{G}}} \widehat{f}(\xi,\zeta) \,\widehat{g}(-\xi-\xi',-\zeta-\zeta') \,\widehat{h}(\xi',\zeta') \,\widehat{\chi}(\xi+\xi'-\zeta') \,\widehat{\chi}(\zeta+\xi'+\zeta'), \label{eq:weiFour2} \\
\Lambda_{\chi}(f,g,h) & = \sum_{\xi,\zeta,\xi',\zeta'\in\widehat{\mathbb{G}}} \widehat{f}(\xi,\zeta) \,\widehat{g}(\xi',\zeta') \,\widehat{h}(-\xi-\xi',-\zeta-\zeta') \,\widehat{\chi}(\zeta-\xi'+\zeta') \,\widehat{\chi}(-\xi-\xi'-\zeta') \label{eq:weiFour3}
\end{align}
\end{subequations}
}
for $f,g,h,\chi$ as before.
Indeed, substituting $x'=x+s$ and $y'=y+t$ we see that $\Lambda_{\chi}(f,g,h)$ is the scalar product of complex-valued functions
\[ (x,y,x',y') \mapsto f(x,y) \,g(x',y') \]
and
\[ (x,y,x',y') \mapsto \overline{h(x+y-y', -x+y+x') \,\chi(-x+x') \,\chi(-y+y')} \]
on the group $\mathbb{G}^4$.
Since the latter function has Fourier transform
\[ \widehat{\mathbb{G}^4}\cong\widehat{\mathbb{G}}^4\to\mathbb{C},\quad (\xi,\zeta,\xi',\zeta') \mapsto \overline{\widehat{h}(-\xi-\xi',-\zeta-\zeta') \, \widehat{\chi}(\zeta-\xi'+\zeta') \,\widehat{\chi}(-\xi-\xi'-\zeta')}, \]
we arrive at \eqref{eq:weiFour3} by applying Plancherel's theorem, i.e., the fact that the Fourier transform is a unitary operator from $\textup{L}^2(\mathbb{G}^4)$ to $\ell^2(\widehat{\mathbb{G}}^4)$.
The same argument guarantees that the right hand side of \eqref{eq:weiFour3} converges absolutely and then the other two formulae, \eqref{eq:weiFour1} and \eqref{eq:weiFour2}, follow by easy changes of variables in the summation over $\widehat{\mathbb{G}}^4$.
A particular case of \eqref{eq:weiFour1} for the constant weight $\chi=\mathbbm{1}$ is worth formulating:
\begin{equation}\label{eq:Four}
\Lambda_{\mathbbm{1}}(f,g,h) = \sum_{\xi,\zeta\in\widehat{\mathbb{G}}} \widehat{f}(-\xi-\zeta,\xi-\zeta) \,\widehat{g}(\xi,\zeta) \,\widehat{h}(\zeta,-\xi).
\end{equation}
Useful inequalities that immediately follow from \eqref{eq:weiFour1}--\eqref{eq:weiFour3} are
{\allowdisplaybreaks
\begin{subequations}
\begin{align}
\big| \Lambda_{\chi}(f,g,h) \big| & \leq \big\|\widehat{f}\big\|_{\ell^{\infty}(\widehat{\mathbb{G}}^2)} \|g\|_{\textup{L}^2(\mathbb{G}^2)} \|h\|_{\textup{L}^2(\mathbb{G}^2)} \big\|\widehat{\chi}\big\|_{\ell^1(\widehat{\mathbb{G}})}^2, \label{eq:Lest2a} \\
\big| \Lambda_{\chi}(f,g,h) \big| & \leq \|f\|_{\textup{L}^2(\mathbb{G}^2)} \big\|\widehat{g}\big\|_{\ell^{\infty}(\widehat{\mathbb{G}}^2)} \|h\|_{\textup{L}^2(\mathbb{G}^2)} \big\|\widehat{\chi}\big\|_{\ell^1(\widehat{\mathbb{G}})}^2, \label{eq:Lest2b} \\
\big| \Lambda_{\chi}(f,g,h) \big| & \leq \|f\|_{\textup{L}^2(\mathbb{G}^2)} \|g\|_{\textup{L}^2(\mathbb{G}^2)} \big\|\widehat{h}\big\|_{\ell^{\infty}(\widehat{\mathbb{G}}^2)} \big\|\widehat{\chi}\big\|_{\ell^1(\widehat{\mathbb{G}})}^2. \label{eq:Lest2c}
\end{align}
\end{subequations}
}
For instance, equality \eqref{eq:weiFour2} and the Cauchy--Schwarz inequality in $\xi,\zeta,\xi',\zeta'$ give
\begin{align*}
\big| \Lambda_{\chi}(f,g,h) \big| \leq \big\|\widehat{g}\big\|_{\ell^{\infty}(\widehat{\mathbb{G}}^2)} & \bigg( \sum_{\xi,\zeta\in\widehat{\mathbb{G}}} \big| \widehat{f}(\xi,\zeta) \big|^2 \Big( \sum_{\xi',\zeta'\in\widehat{\mathbb{G}}} \big|\widehat{\chi}(\xi+\xi'-\zeta')\big| \,\big|\widehat{\chi}(\zeta+\xi'+\zeta')\big| \Big) \bigg)^{1/2} \\
\times & \bigg( \sum_{\xi',\zeta'\in\widehat{\mathbb{G}}} \big| \widehat{h}(\xi',\zeta') \big|^2 \Big( \sum_{\xi,\zeta\in\widehat{\mathbb{G}}} \big|\widehat{\chi}(\xi+\xi'-\zeta')\big| \,\big|\widehat{\chi}(\zeta+\xi'+\zeta')\big| \Big) \bigg)^{1/2}
\end{align*}
and then Plancherel's theorem on $\mathbb{G}^2$ deduces \eqref{eq:Lest2b}.
Here is where we need the $2$-divisibility assumption on $\mathbb{G}$. It is equivalent with the fact that the map $\widehat{\mathbb{G}}\to\widehat{\mathbb{G}}$, $\xi\mapsto 2\xi$ is injective and here we need it to argue that the map
\[ \widehat{\mathbb{G}}^2 \to \widehat{\mathbb{G}}^2,\quad (\xi',\zeta') \mapsto (\xi+\xi'-\zeta',\zeta+\xi'+\zeta') \]
takes each value of $\widehat{\mathbb{G}}^2$ at most once for each fixed choice of $\xi,\zeta\in\widehat{\mathbb{G}}$.
On the other hand, already the definition \eqref{eq:defLambda} gives easy inequalities
{\allowdisplaybreaks
\begin{subequations}
\begin{align}
\big| \Lambda_{\chi}(f,g,h) \big| & \leq \|f\|_{\textup{L}^1(\mathbb{G}^2)} \|g\|_{\textup{L}^{\infty}(\mathbb{G}^2)} \|h\|_{\textup{L}^{\infty}(\mathbb{G}^2)} \|\chi\|_{\textup{L}^1(\mathbb{G})}^2, \label{eq:Lest3a} \\
\big| \Lambda_{\chi}(f,g,h) \big| & \leq \|f\|_{\textup{L}^{\infty}(\mathbb{G}^2)} \|g\|_{\textup{L}^1(\mathbb{G}^2)} \|h\|_{\textup{L}^{\infty}(\mathbb{G}^2)} \|\chi\|_{\textup{L}^1(\mathbb{G})}^2, \label{eq:Lest3b} \\
\big| \Lambda_{\chi}(f,g,h) \big| & \leq \|f\|_{\textup{L}^{\infty}(\mathbb{G}^2)} \|g\|_{\textup{L}^{\infty}(\mathbb{G}^2)} \|h\|_{\textup{L}^1(\mathbb{G}^2)} \|\chi\|_{\textup{L}^1(\mathbb{G})}^2. \label{eq:Lest3c}
\end{align}
\end{subequations}
}
For example,
\[ \big| \Lambda_{\chi}(f,g,h) \big| \leq \|f\|_{\textup{L}^{\infty}(\mathbb{G}^2)}  \|h\|_{\textup{L}^{\infty}(\mathbb{G}^2)}  \int_{\mathbb{G}} \int_{\mathbb{G}} \Big( \int_{\mathbb{G}} \int_{\mathbb{G}} |g(x+s,y+t)| \,\textup{d}\mu(x) \,\textup{d}\mu(y) \Big)  \,|\chi(s)| \,|\chi(t)| \,\textup{d}\mu(s) \,\textup{d}\mu(t), 
 \]
which confirms \eqref{eq:Lest3b}.

We are closely following the elegant ``energy pigeonholing'' proof of Roth's theorem by Tao \cite{T14:Roth}, appearing to some extent in the papers by Bourgain \cite{Bo86:Roth} and Green \cite{G05:3ap}.
Thus, we will be concise, performing only a few necessary modifications.

For $M=2\lfloor 10^4\varepsilon^{-2}\rfloor+4$ we recursively construct an increasing sequence $(S_i)_{i=0}^{M}$ of finite subsets of $\widehat{\mathbb{G}}$ and a decreasing sequence $(\delta_j)_{j=0}^{3M+2}$ of positive numbers.
We define them in the order
\[ S_0, \delta_0, \delta_1, \delta_2, S_1, \delta_3, \delta_4, \delta_5, \ldots, S_M, \delta_{3M}, \delta_{3M+1}, \delta_{3M+2}, \]
i.e., each object in this list depends only on preceding ones.
We also denote $\nu_i:= \nu_{S_i,\delta_{3i}}$.

Let the initial set $S_0$ consist of the nul-character and, if $\mathbb{G}=\mathbb{Z}/N\mathbb{Z}$ is a cyclic group, also the character $\mathbb{G}\to\mathbb{R}/\mathbb{Z}$, $x+N\mathbb{Z}\mapsto x/N+\mathbb{Z}$, which distinguishes the points of $\mathbb{G}$. Consequently, in the latter case we know that each Bohr set $B^{\mathbb{G}}(S_i,\rho)$ will certainly be contained in $\{-\lfloor \rho N\rfloor,\ldots,\lfloor \rho N\rfloor\} + N\mathbb{Z}$. Also set $\delta_0:=\varepsilon/10^4$.
The recurrence relation used to define the sets $S_i$ is the following: let $S_{i+1}\subseteq\widehat{\mathbb{G}}$ be the smallest set such that
\[ S_{i+1} \supseteq S_i \cup \big\{\xi\in\widehat{\mathbb{G}} : \big|\widehat{\nu_i}(\xi)\big|\geq \varepsilon/100 \big\} \]
and
\[ S_{i+1}^2 \supseteq \big\{(\xi,\zeta)\in\widehat{\mathbb{G}}^2 : \big|\widehat{f}(\xi,\zeta)\big|\geq \delta_{3i+2} \big\}. \]
In particular, cardinality of the set $S_i$ will be controlled by $\delta_0,\ldots,\delta_{3i-1}$.
We insist on containing large Fourier coefficients of $f$ in a square $S_{i+1}^2$ being motivated by
\[ B^{\mathbb{G}}(S_i,\rho)^2 \subseteq B^{\mathbb{G}^2}(S_i^2,2\rho) \subseteq B^{\mathbb{G}}(S_i,2\rho)^2 \]
for any $\rho>0$.
The numbers $\delta_j$ will be chosen later, in different ways according to the remainder that $j$ gives when divided by $3$.

By the Plancherel formula and the pigeonhole principle there exists an index $0\leq i\leq M-2$ such that
\begin{equation}\label{eq:proof1}
\sum_{(\xi,\zeta)\in S_{i+2}^2\setminus S_i^2} \big|\widehat{f}(\xi,\zeta)\big|^2 \leq \frac{\varepsilon^2}{10^4}.
\end{equation}
Define
\[ f_0 := f\ast(\nu_i\otimes\nu_i),\quad f_1 := f\ast(\nu_{i+1}\otimes\nu_{i+1}-\nu_i\otimes\nu_i),\quad f_2 := f - f\ast(\nu_{i+1}\otimes\nu_{i+1}) \]
and observe that $f_0$ still takes values in $[0,1]$, while $f_1$ and $f_2$ take values in $[-1,1]$.
By
\[ \big|\widehat{f_2}(\xi,\zeta)\big| \leq \big|\widehat{f}(\xi,\zeta)\big| \int_{\mathbb{G}} \int_{\mathbb{G}} \big|1 - e^{-2\pi \mathbbm{i} \xi(x)}e^{-2\pi \mathbbm{i}\zeta(y)}\big| \,\nu_{i+1}(x) \,\nu_{i+1}(y) \,\textup{d}\mu(x) \,\textup{d}\mu(y) \]
we see that $|\widehat{f_2}(\xi,\zeta)|$ can be estimated either by the definition of $S_{i+1}$ for pairs $(\xi,\zeta)\not\in S_{i+1}^2$, or by
the fact that $\nu_{i+1}$ vanishes outside of $B^{\mathbb{G}}(S_{i+1},\delta_{3i+3})$ and
$|1-e^{-2\pi \mathbbm{i} \xi(x)}|,|1-e^{-2\pi \mathbbm{i} \zeta(y)}|\leq 2\pi\delta_{3i+3}$
for $x,y\in B^{\mathbb{G}}(S_{i+1},\delta_{3i+3})$ and pairs $(\xi,\zeta)\in S_{i+1}^2$.
Thus, $\delta_{3i+3}$ can be chosen sufficiently small so that
\begin{equation}\label{eq:proof2}
\big\|\widehat{f_2}\big\|_{\ell^{\infty}(\widehat{\mathbb{G}}^2)} \leq 2\delta_{3i+2}.
\end{equation}
Next, we write
\[ \|f_1\|_{\textup{L}^2(\mathbb{G}^2)}^2 = \sum_{(\xi,\zeta)\in\widehat{\mathbb{G}}^2} \big|\widehat{f}(\xi,\zeta)\big|^2 \big| \widehat{\nu_{i+1}}(\xi) \widehat{\nu_{i+1}}(\zeta) - \widehat{\nu_{i}}(\xi) \widehat{\nu_{i}}(\zeta) \big|^2 \]
and split the sum according to whether $(\xi,\zeta)\in S_{i+2}^2\setminus S_i^2$, or $(\xi,\zeta)\not\in S_{i+2}^2$, or $(\xi,\zeta)\in S_i^2$.
The first part of the sum is controlled using \eqref{eq:proof1}, the second one by $|\widehat{\nu_{i}}(\xi)|,|\widehat{\nu_{i+1}}(\xi)|\leq\varepsilon/100$ or $|\widehat{\nu_{i}}(\zeta)|,|\widehat{\nu_{i+1}}(\zeta)|\leq\varepsilon/100$, following from the definition of $S_{i+2}$, while for the third one we expand the definition of $\widehat{\nu_{i}}$ and $\widehat{\nu_{i+1}}$ as before. We conclude
\begin{equation}\label{eq:proof3}
\|f_1\|_{\textup{L}^1(\mathbb{G}^2)} \leq \|f_1\|_{\textup{L}^2(\mathbb{G}^2)} \leq \frac{\varepsilon}{27}.
\end{equation}
From
\[ |f_0(x+s,y+t) - f_0(x,y)| \leq \int_{\mathbb{G}} \int_{\mathbb{G}} |f(x-u,y-v)| \,\big| \nu_{i}(u+s) \nu_{i}(v+t) - \nu_{i}(u) \nu_{i}(v) \big| \,\textup{d}\mu(u) \,\textup{d}\mu(v) \]
and
\[ \big| \nu_{i}(u+s) \nu_{i}(v+t) - \nu_{i}(u) \nu_{i}(v) \big| \leq \frac{\|\nu_i\|_{\textup{L}^{\infty}(\mathbb{G})}^2}{\delta_{3i}}
\Big( \max_{\xi\in S_i} \| \xi(s) \|_{\mathbb{R}/\mathbb{Z}} + \max_{\zeta\in S_i} \| \zeta(t) \|_{\mathbb{R}/\mathbb{Z}} \Big) \]
we also see that $\delta_{3i+1}$ can be chosen sufficiently small so that
\begin{equation}\label{eq:proof4}
\sup_{\substack{x,y\in\mathbb{G}\\ s,t\in B^{\mathbb{G}}(S_i,2\delta_{3i+1})}}|f_0(x+s,y+t) - f_0(x,y)| \leq \frac{\varepsilon}{54}.
\end{equation}

Finally, we take $S:= S_i$ and $\rho:= \delta_{3i+1}$.
Using multilinearity we split
\[ \Lambda_{\chi_{S,\rho}}(f,f,f) = \Lambda_{\chi_{S,\rho}}(f_0+f_1+f_2,f_0+f_1+f_2,f_0+f_1+f_2) \]
into $27$ terms.
The terms involving $f_1$ are estimated by $\varepsilon/27$ using \eqref{eq:Lest3a}--\eqref{eq:Lest3c} and \eqref{eq:proof3}.
Among the remaining terms, those that contain $f_2$ are bounded using \eqref{eq:Lest2a}--\eqref{eq:Lest2c} and \eqref{eq:proof2} by
\[ 2\delta_{3i+2} \big\|\widehat{\chi_{S,\rho}}\big\|_{\ell^1(\widehat{\mathbb{G}})}^2
= 2\delta_{3i+2} \big\|\widehat{\nu_{S,\rho}}\big\|_{\ell^2(\widehat{\mathbb{G}})}^4
= 2\delta_{3i+2} \|\nu_{S_i,\delta_{3i+1}}\|_{\textup{L}^2(\mathbb{G})}^4, \]
which is at most $\varepsilon/27$ provided that $\delta_{3i+2}$ is sufficiently small in terms of $|S_{i}|$ and $\delta_{3i+1}$.
It remains to observe that $\Lambda_{\chi_{S,\rho}}(f_0,f_0,f_0)$ differs by at most $\varepsilon/27$ from
\[ \int_{\mathbb{G}} \int_{\mathbb{G}} \int_{\mathbb{G}} \int_{\mathbb{G}} f_0(x,y)^3 \,\chi_{S,\rho}(s) \,\chi_{S,\rho}(t) \,\textup{d}\mu(s) \,\textup{d}\mu(t) \,\textup{d}\mu(x) \,\textup{d}\mu(y)
= \int_{\mathbb{G}} \int_{\mathbb{G}} f_0(x,y)^3 \,\textup{d}\mu(x) \,\textup{d}\mu(y), \]
thanks to \eqref{eq:proof4}.
The last display is greater than or equal to the third power of
\[ \int_{\mathbb{G}} \int_{\mathbb{G}} \big( f\ast(\nu_i\otimes\nu_i) \big)(x,y) \,\textup{d}\mu(x) \,\textup{d}\mu(y) = \int_{\mathbb{G}} \int_{\mathbb{G}} f(x,y) \,\textup{d}\mu(x) \,\textup{d}\mu(y). \]
This completes the proof of Theorem~\ref{thm:compact}.


\section{Proof of Theorem~\ref{thm:expsmall}}
\label{sec:expsmall}
This time we are closely following the approach by Tao from an expository note \cite{T16:poly}.

Observe that the configuration \eqref{eq:triangle} can equivalently be written as
\[ (x,y), \ (x',y'), \ (x'',y'') \]
with constraints $x-y-x'+y''=\mathbf{0}$ and $x+y-y'-x''=\mathbf{0}$ in $\mathbb{F}_q^n$.
For a set $A$ as in the theorem formulation we have the equality
\begin{equation}\label{eq:Fidentity}
\delta_{\mathbf{0}}(x-y-x'+y'') \delta_{\mathbf{0}}(x+y-y'-x'') = \sum_{(a,b)\in A} \delta_{(a,b)}(x,y) \delta_{(a,b)}(x',y') \delta_{(a,b)}(x'',y'')
\end{equation}
of functions $A^3\to\mathbb{F}_q$ in variables $(x,y),(x',y'),(x'',y'')$.
Just as in \cite{T16:poly} we define the \emph{(slice) rank} of a function $F\colon A^3\to\mathbb{F}_q$ to be the minimal number of terms of the form
\[ f(x,y)g(x',y',x'',y'') \quad\text{or}\quad f(x',y')g(x,y,x'',y'') \quad\text{or}\quad f(x'',y'')g(x,y,x',y') \]
that sum up to $F((x,y),(x',y'),(x'',y''))$.
From \cite[Lemma~1]{T16:poly} we know that the right hand side of \eqref{eq:Fidentity} has slice rank precisely $|A|$.

On the other hand, if $x=(x_i)_{i=1}^{n}$, $x'=(x'_i)_{i=1}^{n}$, etc., then the left hand side of \eqref{eq:Fidentity} can be written as
\begin{equation}\label{eq:Fidentity2}
\prod_{i=1}^n \big(1-(x_i-y_i-x'_i+y''_i)^{q-1}\big) \big(1-(x_i+y_i-y'_i-x''_i)^{q-1}\big).
\end{equation}
Note that the polynomial \eqref{eq:Fidentity2} is a linear combination of the terms
\[ \prod_{i=1}^n (x_i-y_i)^{\alpha_i} (x_i+y_i)^{\beta_i} (x'_i)^{\alpha'_i} (y'_i)^{\beta'_i} (y''_i)^{\alpha''_i} (x''_i)^{\beta''_i} \]
for some $\alpha_i,\beta_i,\alpha'_i,\beta'_i,\alpha''_i,\beta''_i\in\{0,1,\ldots,q-1\}$ such that
\[ \sum_{i=1}^n (\alpha_i+\beta_i+\alpha'_i+\beta'_i+\alpha''_i+\beta''_i) \leq 2(q-1)n. \]
We divide those terms into three groups, one of them determined by tuples of exponents such that
\begin{equation}\label{eq:firstgroup}
\sum_{i=1}^n (\alpha_i+\beta_i) \leq \frac{2(q-1)n}{3}
\end{equation}
and the other two groups defined analogously. The terms in the first group can be organized into a sum of $D$ products of the form
\[ \Big(\prod_{i=1}^n (x_i-y_i)^{\alpha_i} (x_i+y_i)^{\beta_i}\Big) g(x',y',x'',y''), \]
where $D$ is the number of tuples $(\alpha_1,\ldots,\alpha_n,\beta_1,\ldots,\beta_n)\in\{0,1,\ldots,q-1\}^{2n}$ satisfying \eqref{eq:firstgroup}.
Note that $\mathop\textup{char}\mathbb{F}_q=2$ is allowed, as this degeneracy can only decrease the count.
Using a standard trick of estimating a restricted sum by a weighted sum, we can bound $D$ for any $0<t\leq 1$ as
\[ D \leq \sum_{\alpha_1,\ldots,\alpha_n,\beta_1,\ldots,\beta_n=0}^{q-1} t^{\sum_{i=1}^n (\alpha_i+\beta_i) - \frac{2(q-1)n}{3}} = \varphi_q(t)^{2n}, \]
where $\varphi_q(t) := (1+t+\cdots+t^{q-1})t^{-(q-1)/3}$. Thus, the slice rank of the left hand side of \eqref{eq:Fidentity} is at most
\[ 3D \leq 3\Big(\min_{t\in(0,1]}\varphi_q(t)\Big)^{2n}. \]
The above minimum is certainly less than $q$, because $\varphi_q(1)=q$ and $\varphi'_q(1)=q(q-1)/6>0$.
This completes the proof of Theorem~\ref{thm:expsmall}.


\section*{Acknowledgments}
This work is supported in part by the \emph{Croatian Science Foundation} project UIP-2017-05-4129 (MUNHANAP).
The author is grateful to Aleksandar Bulj and Rudi Mrazovi\'{c} for several comments and useful discussions.
The author is also thankful for numerous comments and useful suggestions by an anonymous referee.


\end{document}